\documentclass[fleqn]{article}

\usepackage{makeidx}         
\usepackage{graphicx}        
\usepackage{multicol}        
\usepackage[bottom]{footmisc}
\usepackage{amsfonts}
\usepackage{epsfig}
\usepackage{stmaryrd}

\author{Francisco Bernal$^{1}$, Gail Guti\'errez$^{2}$ \\
\ $^{1}$Universidad Carlos III, 28911 Legan\'es (Madrid) \\
$^{2}$Universidad Pontificia Bolivariana (Medell\'in)}
\title{SOLVING DELAY DIFFERENTIAL EQUATIONS THROUGH RBF COLLOCATION}
\begin{document}
\maketitle

\begin{abstract}

A general and easy-to-code numerical method based on radial basis functions (RBFs)
collocation is proposed for the solution of delay differential equations (DDEs).
It relies on the interpolation properties of infinitely smooth RBFs, which allow
for a large accuracy over a scattered and relatively small discretization support.
Hardy's multiquadric is chosen as RBF and combined with the Residual Subsampling
Algorithm of Driscoll and Heryudono for support adaptivity. The performance of the method 
is very satisfactory, as demonstrated over a cross-section of benchmark DDEs, and by 
comparison with existing general-purpose and specialized numerical schemes for DDEs.     
\end{abstract}

\section{Introduction}

In this work, we present a general numerical approach for solving DDEs based on
the RBF collocation method invented by Kansa \cite{Kansa90a}\cite{Kansa90b}, also
known as Kansa's method. Due to its many advantages (which include
superior interpolation accuracy, spectral convergence, robustness with respect
to the discretization support, and ease of coding), Kansa's method is
becoming increasingly popular for the solution of ordinary and partial differential 
equations (ODEs and PDEs, respectively). Its performance in the solution of
DDEs, however, has scarcely been explored, with the exception of 
a recent paper on the solution
of neutral DDEs with multiquadrics \cite{Karimi08}.  
This paper is organized as follows. In Section 2, Kansa's
method is adapted to a general formulation of (first order) DDEs. The
basic algorithm is further improved by the inclusion of several heuristic
observations concerning the tunable shape parameter which appears in the
multiquadric RBF, and by the residual subsampling algorithm (RSA) by
Driscoll and Heryuodono \cite{Driscoll07}. The RSA is at the core of the
high accuracy attained by the multiquadrics interpolant. Section 2 is closed by 
some remarks concerning the solution of nonlinear problems with
Kansa's method.
Section 3 tests the proposed method against a cross-section
of benchmark problems taken from the literature. As we shall see, not
only does Kansa's method attain excellent results in well-understood (first order)
DDEs, but also in the less explored neutral and higher-order DDEs 
-which may offer an additional tool for looking into this kind of problems.
Finally, Section 4 concludes the paper.     

\section{Solving linear DDEs through Kansa's method}

Consider the following linear DDE

\begin{equation}
y'(x) - p(x)y(x) - q(x)y[x-\tau(x)] = s(x) \, \textrm{ if } x\in[a,b]
\end{equation}
\begin{equation}
y(x)= h(x)  \, \textrm{ if } x\leq a
\label{F:general}
\end{equation}

It will be convenient to split (\ref{F:general}) into a DDE and an ODE

\begin{equation}
y'(x) - p(x)y(x) - q(x)y[x-\tau(x)] = s(x) \qquad \textrm{ if } x-\tau(x)>a
\label{F:splitdde1}
\end{equation}
\begin{equation}
y'(x) - p(x)y(x) = q(x)h[x-\tau(x)] + s(x) \qquad \textrm{ if } x-\tau(x)<a
\label{F:splitdde2}
\end{equation}
\begin{equation}
y(a) = h(a)
\label{F:splitdde3}
\end{equation}

Discretize $[a,b]$ into a set $N$ scattered nodes $\xi= \{x_j, j=1...N\}$ (with
$x_1=a$ and $x_N=b$), and
consider as well the outside point $x_0= a - \lambda,\,\lambda>0$. 
We seek an approximate solution to (\ref{F:splitdde1})-(\ref{F:splitdde3}) in the form of an expansion
of $N+1$ RBFs $\phi_j(r)$:
\begin{equation}
y(x)= \sum_{j=0}^{j=N} \alpha_j \phi(\parallel x-x_j \parallel)  
\label{F:RBFinterpolant}
\end{equation} 

The addition of an RBF at $x_0$ allows to enforce \textit{both} the
initial condition \textit{and} the DDE at $x=a$, thus contributing to
the accuracy (this is the PDECB strategy discussed in \cite{Fedoseyev00}). Once the coefficients $\alpha_j$
are available, the approximate RBF solution can be reconstructed 
anywhere in $[a,b]$. In order to solve for the coefficients, 
(\ref{F:splitdde1})-(\ref{F:splitdde3}) are enforced over (\ref{F:RBFinterpolant}) on a set of 
collocation $N$ nodes, usually $\xi$.
Notice that no equation is collocated on $x_0$, but two of them
are on $x_1=a$. For $i=1,\ldots,N$, this leads to the linear system of dimension $N+1$

\begin{equation}
\sum_{j=0}^{j=N} \{ \phi '_j(r_{ij}) - p(x_i)\phi_j(r_{ij}) - 
q(x_i)\phi_j(||x_i-\tau(x_i)-x_j||) \}=
s(x_i)\textrm{  if } x_i-\tau(x_i) > a
\label{F:colsys1}
\end{equation}
\begin{equation}
\sum_{j=0}^{j=N} \{ \phi '_j(r_{ij}) - p(x_i)\phi_j(r_{ij})\} \alpha_j= 
q(x_i)h[x_i-\tau(x_i)] + s(x_i)
\textrm{  if } x_i-\tau(x_i) \le a 
\label{F:colsys2}
\end{equation}
\begin{equation}
\sum_{j=0}^{j=N} \alpha_j \phi_j(r_{ij})= h(a)
\textrm{  if } x_i=a
\label{F:colsys3}
\end{equation}

where $r_{ij}= || x_i - x_j ||$. In the remainder of this paper we will restrict 
ourselves to the well-tested Hardy's multiquadric (MQ),
\begin{equation}
\phi_j(r_j)= \sqrt{\parallel x-x_j \parallel^2 + c_j^2}
\label{F:MQ}
\end{equation}
whose derivative is
\begin{equation}
\phi '_j(r_j)= \frac{x-x_j}{\sqrt{\parallel x-x_j \parallel^2 + c_j^2}}
\end{equation}
as the RBF of choice. The shape of the MQ depends on the free parameter $c_j$ (hence the name of \textit{shape parameter} for it). The fact that the MQ has global support leads to fully populated matrices.
It is a hallmark of Kansa's method that the best accuracy can only be obtained at
the expense of extreme ill-conditioning, as will be discussed next. In order
to improve stability, the direct inversion of the linear system (\ref{F:colsys1})-(\ref{F:colsys3})
has been replaced by the use of Penrose's pseudoinverse. 

\subsection{Choosing the shape parameters $c_j$}

Although the accuracy of the interpolant (\ref{F:RBFinterpolant}) is largely influenced by the 
values $c_j,\,j=0,\ldots,N+1$,
theoretical results regarding the choice of an 'optimal' set of values are not yet available,
and heuristic rules must be used instead, which mostly address the homogeneous case $c_j=c$. In this case, the convergence
rate of the error of the interpolant (\ref{F:RBFinterpolant}) has been proven to go as $\Lambda ^{c/h}$ in interpolation
problems \cite{Madych92}, and has been shown to obey $\Lambda ^{\sqrt{c}/h}$ in elliptic PDEs \cite{Cheng03},
where $0<\Lambda<1$ and $h$ is the distance between nodes. Therefore, the
accuracy could seemingly be improved at no computational cost by increasing $c$. However, as 
$c \shortrightarrow \infty$, the MQ profile becomes increasingly flatter and the collocation system
(\ref{F:colsys1})-(\ref{F:colsys3}) becomes extremely ill-conditioned, dictating in practice a limit for the accuracy attainable at
a given resolution $h$ and machine precision. A trade-off principle arises between accuracy and stability, 
which is actually common to all parameter-dependent RBFs, not only MQs \cite{Schaback95}. Optimal results are obtained by pushing 
$c$ as large as
possible before incurring in numerical instability. Since in the approximation of differential equations the
exact solution is unknown, other estimators are used instead, which replicate the behavior of the error curves with $c$
and are available in run time. Examples are the 'leave-one-out' strategy 
\cite{Rippa99}\cite{Fasshauer07c}, or the residual to the ODE/PDE \cite{Cheng03}. For instance, in (\ref{F:splitdde1})-(\ref{F:splitdde3}),the pointwise residual is defined as 
\begin{equation}
R(x)= s(x) - \sum_{j=0}^{j=N}\alpha_j\phi '_j(x) + 
p(x)\sum_{j=0}^{j=N}\alpha_j\phi_j(x) + q(x)\sum_{j=0}^{j=N}\alpha_j\phi '_j[x-\tau(x)] 
\end{equation}
The case where $c$ is center-dependent has been less 
investigated, although it may outperform MQ collocation with constant $c$, as shown in a 
numerical investigation by Kansa and Carlson \cite{Kansa92}. Carlson and Foley showed that $c_j$ is related to the curvature of the
function to be interpolated at $x\approx x_j$ \cite{Carlson91}. In \cite{Hon97}, Hon and Mao let $c_j=Mj+b$, 
where $j$ is the center index and $M$ and $b$ are chosen so that the condition number $\kappa$ is about $10^{16}$. 
In \cite{Wertz06},  Wertz {\em et al.} reported improved accuracy in a 2D problem if 
$c_j\gg c_0$ for $(x_j,y_j)\in \partial \Omega$ and $c_j= \mu(1+\gamma(-1)^j)$ if $(x_j,y_j) \in \Omega$, for some 
constants $\mu$ and $\gamma$. These findings were confirmed in the 1D case in a later work by 
Fornberg and Zuev \cite{Fornberg07}. Another common strategy has been to set $c$ proportional to the distance to the 
closest node in the point set, {\em v.g.} in \cite{Driscoll07}.

\subsection{Extension to nonlinear DDEs}           
 
In the case that the DDE is nonlinear, or that the lagged argument is a function of
the solution itself (\textit{a state-delay DDE}), the collocation of the interpolant
(\ref{F:RBFinterpolant}) leads to a system of nonlinear algebraic equations for
the unknowns $\alpha_0,\ldots,\alpha_N$. Let us write this system as
\begin{eqnarray}
F_0(\alpha_0,\ldots,\alpha_N)= 0 \nonumber\\
\vdots \\
F_N(\alpha_0,\ldots,\alpha_N)= 0 \nonumber
\end{eqnarray}
In order to solve $\vec{F}=\vec{0}$, a gradient-based method may be used. In the
MATLAB routine \textit{fsolve}, the user can choose between providing the analytical
Jacobian $J$ to the solver,
\begin{equation}
J=
\left(
\begin{array}{ccc}
\frac{\partial F_0}{\partial \alpha_0} & \ldots & \frac{\partial F_0}{\partial \alpha_N} \\
\vdots & \ddots & \vdots \\
\frac{\partial F_N}{\partial \alpha_0} & \ldots & \frac{\partial F_N}{\partial \alpha_N} \\
\end{array}
\right)
\end{equation}
or allowing it to construct $J$ based on finite differences. In order to keep the
implementation of Kansa's method as simple and general as possible, we
have only explored the latter possibility. However, there is a practical
drawback: while it is well known
that the convergence of Newton-type methods is very sensitive to the condition
number of the Jacobian, RBF interpolation needs to push $\kappa$ for the
best accuracy, often beyond the ill-condition threshold (which is $\kappa\approx 10^{14}$
in our MATLAB environment). Consequently, we have used instead a trust-region
method (that of Powell's \cite{Powell70}) in our numerical experiments, with good results.
Nevertheless, the condition number must be kept lower than in linear
DDEs in order to guarantee convergence, which is likely to prevent optimal
accuracy as well.

\subsection{Adaptive selection of nodes}
Another important yet open issue in Kansa's method is the optimal number and location of RBF 
centers/collocation nodes. We will restrict ourselves to the case where both point sets are identical 
(save for the extra RBF center at $x_0$ added in order to enforce the equation at $x=a$). In 1D problems, 
regular grids are often preferred for simplicity, although experimental evidence suggests that the optimal placement of nodes is problem-dependent, {\em i.e.} is determined by the function to be interpolated. We will be using an algorithm introduced by Driscoll and Heryuodono \cite{Driscoll07} which works well in practice, both for interpolation and differential equations and not only in 1D. The idea is to monitor the residual $R$ to the differential equation at midpoints and iteratively refine the point set until $R$ drops below some user-defined threshold. The reader is referred to the original paper for details. Here, we present a slightly modified version of the algorithm which we have preferred.
\\
\\
\textit{Residual Subsampling Algorithm (RSA)}  \\ 
\begin{itemize}
\item Initially, discretize  $[a,b]$ into a grid of $N^{(0)}$ nodes with spacing $\Delta=(b-a)/(N^{(0)}-1)$.
Define $x_j=a+(j-1)\Delta\},\,j=1,\ldots,N^{(0)}$, $\xi^{(0)}= \{x_j\}$, and $x_0=a-\Delta$. The $N^{(0)}+1$ starting MQ centers are the set 
$x_0 \cup \xi^{(0)}$. Define the values of the adjustable parameters $\lambda>0$, $\mu>0$, $\gamma>0$, $\eta>0$, 
$\theta_{max}>\theta_{min}>0$, and \textit{itmax}. 
\item For $k=0,\ldots$ until $\max |R_j^{(k)}|<\theta_{max}$ or $k>$\textit{itmax}  
\begin{itemize}
\item Distribute the shape parameters as $c_0 = c_{N^{(k)}}= \lambda \mu d_1$, and $c_j= \mu d_j[1+\gamma(-1)^j], j=1,\ldots,N^{(k)}-1$, where $d_j$ is the distance to the closest collocation node from $x_j$.
\item Compute set of midpoints $z_j=(x_j+x_{j+1})/2, j=1...N^{(k)}-1$.
\item Solve the DDE through Kansa's method with $y(x)= \displaystyle{\sum_{j=0}^{j=N^{(k)}}}\alpha_j \phi_j(||x-x_j||)$.
\item Compute the residuals $\{R_j^{(k)}\}$ to the differential equation at midpoints.
\item Set $\Theta^{(k)}= \max\big(\theta_{max},\displaystyle{\max_{j=1...N^{(k)}-1}}|R_j^{(k)}/\eta|\big)$.
\item Define point set $\Xi^{(k)}= \xi^{(k)} \cup \{z_j \textrm{ such that } |R_j^{(k)}|>\Theta^{(k)}\}$.
\item Delete points $x_i, i=2...N^{(k)}-1$ such that $|R_{i-1}^{(k)}|<\theta_{min}>|R_{i+1}^{(k)}|$ from $\Xi{^(k)}$.
\item Let $\xi^{(k+1)}=\Xi^{(k)}$ and $\{x_1,\ldots,x_{N^{(k+1)}}\}=\xi^{(k+1)}$.
\item Update $\{d_j\}$ for $j=1...N^{(k+1)}$ 
\item Consider the set of $1+N^{(k+1)}$ MQs centered at $x_0 \cup \xi^{(k+1)}$ and iterate.
\end{itemize}
\end{itemize} 

In the above algorithm, the shape parameters are adjusted after each iteration in order
to prevent the condition number from skyrocketing. As further nodes are included, however,
the onset of instability will be eventually reached and the accuracy of the MQ approximation
begins to deteriorate. The only tweakings to the  original RSA in \cite{Driscoll07} are: PDECB, 
the use of the recipe in \cite{Wertz06} in the distribution of $c$'s, and the subtitution of $\theta_{max}$
by $\Theta^{(k)}$ on enlargement of the point set.  

\section{Numerical Examples}
In the remainder of the paper, we will refer to the method described in section 2 as MQCM (multiquadric collocation method). The MQCM is coded in MATLAB 7 running on a laptop with 1.8 GHz CPU and 1 GB RAM.
In this section, the MQCM is tested against a cross-section of
benchmark DDEs taken from the literature. The performance of the MQCM is compared
with that of MATLAB built-in general-purpose routines DDE23 \cite{Shampine01} by Shampine and Thompson, or DDESD \cite{Shampine05} by 
Shampine, which are both based on
Runge-Kutta-type schemes. DDE23 is restricted to
constant delays, while the more recent DDESD can handle variable- and state-delay
equations as well. For Examples 4 and 5, DDENSD has been
used instead of DDESD, which is a routine based on DDESD for DDEs of neutral type. 
The fact that these three programs are written in MATLAB allows for a direct
comparison of error estimates and CPU times with MQCM. In particular, the root mean
squared error is defined as
\begin{equation}
RMS(\epsilon)= \sqrt{\frac{\sum_{i=1}^{i=N_{ev}} [u_{NUM}(z_i)-u_{EX}(z_i)]^2}{N_{ev}}}
\end{equation}  
where $u_{EX}$ is the exact solution, $u_{NUM}$ the approximation yielded by the considered
numerical scheme, $\epsilon$ is the point-wise error, and $z_i,\,i=1,\ldots,N_{ev}=103$ is a set of equispaced evaluation points in
$[a,b]$.
In some of the examples presented, published results of some specialized method for the kind of
DDE considered have been included as further reference. In such cases, not all the
estimators are available for comparison. CPU times, in particular, cannot be directly compared
 -which is denoted by adding an * to the corresponding entry.   
     
In all of the numerical examples which follow, the working parameters for the RSA have been set to
\begin{equation}
\lambda=10, \,\,\,\, \mu=\sqrt{40/N^{(0)}},\,\,\,\, \gamma=0.1,\,\,\,\, \eta= 10,\,\,\,\, 
\theta_{max}=10^{-13},\,\,\,\, \theta_{min}=10^{-14}
\end{equation}
except in Example 5 where $\mu=\sqrt{25/N^{(0)}}$. The initial discretization is
$N^{(0)}=6$ in Examples 1-4, and $N^{(0)}=10$ in Examples 5 and 6.

\subsection{Example 1: Stiff DDE}
\begin{figure}[htb!]
\begin{center}
\includegraphics[height=7cm]{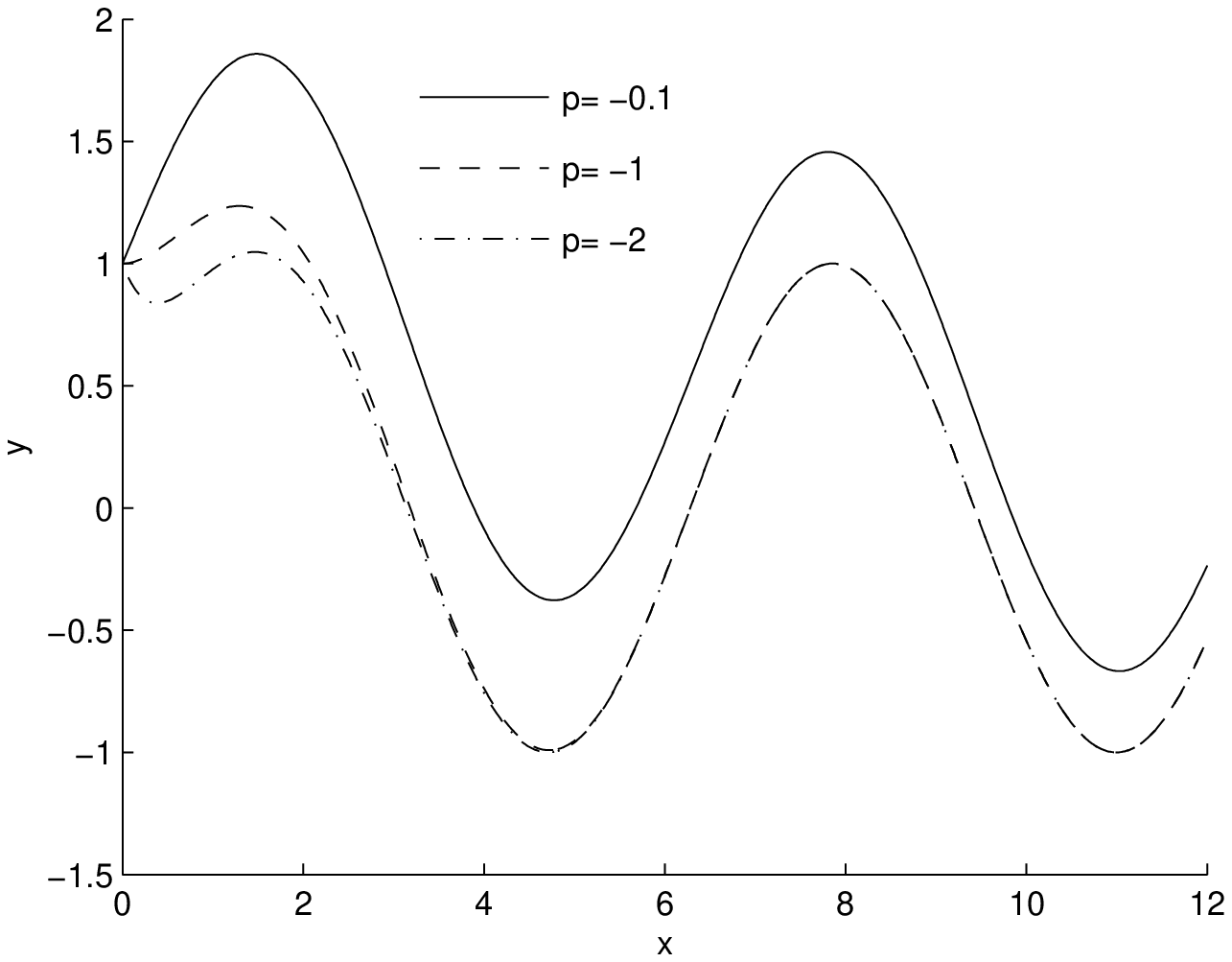}
\caption{Plots of the exact solution of Example 1}
\label{valores_p} 
\end{center}
\end{figure}

Consider the following DDE with a stiffness parameter $p$ (Example 1 in \cite{Ito91}).  

\begin{equation}
\left\{
\begin{array}{l}
y'(x)= Ay(x) + y \big( x-\frac{3\pi}{2} \big) - A\sin (x), \qquad x\in [0,13] \\
y(x)= e^{px} + \sin (x), \qquad x\in \big[ -\frac{3\pi}{2}, 0 \big] \\
\end{array}
\right.
\end{equation}

where $A=p-e^{-3\pi p/2}$. The exact solution is given by $y_{EX}(x)= e^{px}+\sin (x)$.
For $p<0$, the solution consists of a short transient of exponential decay, followed by periodic sinusoidal oscillations (see Fig. \ref{valores_p}). Since the parameter $p$ also enters the equation exponentially, its effect on the stiffness of the problem is dramatic. Table \ref{test_prob1} compares the performance of the MQCM 
with that of DDE23 and with that of the spectral method in
\cite{Ito91} (SPC). In Table \ref{test_prob1}, an entry like $5.1(-15)$ means $5.1\times 10^{-15}$, and so on.  
DoF (degrees of freedom) stands for the size of support of the given discretization 
scheme -the number of MQ centers in the MQCM. The listed results for DDE23 are the best within a reasonable computing time and/or memory restrictions.

The MQCM is barely affected, if anything, by the increasing stiffness of the problem. In fact, the advantages of the MQCM in dealing with stiff ODEs were already reported in \cite{Hon97}. In terms of efficiency, the MQMC outperforms DDE23. The inversion of full matrices required by the MQCM is made up for by the gain in the size of the discretization support. 

On the other hand, the accuracy of the SPC can be improved by increasing the order of the scheme,
as happens in Table \ref{test_prob1} for different $p$. The SPC is more efficient
than the MQCM, but is affected by the increasing value of $p$ (see discussion in \cite{Ito91}). Moreover, it is restricted to constant delays.

\begin{table}[!hbp]	
\centering
\caption{Comparison to other methods (Example 1)}
\label{test_prob1}   
\begin{footnotesize}
\begin{tabular}{|l|l|lll|}
\hline
$p$	&	$x$ & $\epsilon_{MQ}$ & $\epsilon_{DDE23}$ & $\epsilon_{SPC}$ \\
\hline

	&	$3\pi / 4$	&	5.1(-15)	&	1.9(-12)	&	2.6(-6)	\\
	&	$3\pi / 2$	&	6.2(-14)	&	1.5(-12)	&	7-9(-8) \\
-0.1&	$9\pi / 4$	&	9.7(-14)	&	1.4(-12)	&	1.0(-5)	\\
	&	$3\pi$		&	3.5(-14)	&	8.2(-13)	&	3.1(-7) \\
	&	$15\pi / 4$	&	1.6(-13)	&	8.6(-13)	&	8.4(-7) \\

\hline
\multicolumn{2}{|c|}{DoF}  &	261	&	71072	&	27 \\
\hline
\multicolumn{2}{|c|}{RMS($\epsilon$)} & 9.4(-14)	&	2.3(-12) & \\
\hline
\multicolumn{2}{|c|}{CPU} & 16	&	891 & 0.009 *\\
\hline

	&	$3\pi /4$	&	1.9(-13)	&	2.5(-11)	& 	8.3(-9)	\\
	&	$3\pi /2$	&	7.0(-14)	&	1.3(-11)	& 	7.6(-7)	\\
-1	&	$9\pi /4$	&	4.2(-14)	&	6.6(-12)	&	1.5(-8)	\\
	&	$3\pi$		&	6.6(-14)	&	2.0(-14)	&	4.2(-7)	\\
	&	$15\pi /4$	&	6.4(-14)	&	1.9(-11)	&	2.0(-7)	\\

\hline
\multicolumn{2}{|c|}{DoF}  &	254		&	82258	&	33  \\
\hline
\multicolumn{2}{|c|}{RMS($\epsilon$)} & 6.0(-14)	&	3.9(-11) & \\
\hline
\multicolumn{2}{|c|}{CPU} & 26	&	1358 & 0.017 * \\
\hline

	&	$3\pi/4$		&	9.3(-14)	&	2.0(-10)	&	1.3(-10) \\
	&	$3\pi/2$		&	1.4(-13)	&	1.2(-10)	&	1.1(-9)	 \\
-2	&	$9\pi/4$		&	4.2(-14)	&	2.0(-10)	&	2.1(-10) \\
	&	$3\pi$			&	1.8(-14)	&	1.0(-12)	&	1.1(-9)	 \\
	&	$15\pi/4$		&	8.1(-14)	&	1.0(-10)	&	2.1(-10) \\

\hline
\multicolumn{2}{|c|}{DoF}  &	281		&	151122	&	51 \\
\hline
\multicolumn{2}{|c|}{RMS($\epsilon$)} & 1.4(-13)	&	2.1(-10) & \\
\hline
\multicolumn{2}{|c|}{CPU} & 19	 &	5850 &  0.036 *\\
\hline

\end{tabular}
\end{footnotesize}
\end{table}

Table \ref{RSA_prob1} shows the performance of the RSA throughout the iterations for this problem. While it converges on average, the scheme is clearly not monotone. It is surprising that the convergence can be sustained at so high condition numbers. While in our implementation of Kansa's method the ill-conditioning problem is -at least partially- ameliorated by the use of the pseudoinverse (instead of the direct inversion of the matrix), this phenomenon of high accuracy at very high condition numbers has already been reported when smooth functions are interpolated with MQs \cite{Hon97}.

\begin{table}[!hbp]	
\centering
\caption{RSA iterations (Example 1)}
\label{RSA_prob1}   
\begin{footnotesize}
\begin{tabular}{|l|lll|lll|lll|}
\hline
 & \multicolumn{3}{|c|}{$p=-0.1$} & \multicolumn{3}{|c|}{$p=-1$} & \multicolumn{3}{|c|}{$p=-2$} \\
\hline
it & DoF & RMS($\epsilon$) & Condition $\sharp$ & DoF & RMS($\epsilon$) & Condition $\sharp$ & DoF & RMS($\epsilon$) & Condition $\sharp$ \\
\hline
0	&	7	&	0.77 		&	3.3(+10)		&	7	&	0.94	 	&	1.9(+14)		&	7	&	1.04		&	3.9(+17)	\\
1	&	12	&	0.0040	 	&	3.2(+14)		&	10	&	0.0112	 	&	3.4(+14)		&	11	&	0.02		&	1.2(+17)	\\
2	&	15	&	1.6(-6)		&	4.3(+13)		&	14	&	4.0(-5)		&	5.3(+16)		&	13	&	1.6(-4)	 	&	7.3(+17)	\\
3	&	22	&	7.5(-8)		&	1.7(+16)		&	23	&	1.0(-7)		&	4.3(+17)		&	15	&	3.0(-4) 	&	2.3(+17)	\\
4	&	27	&	3.5(-7)		&	4.1(+18)		&	27	&	6.2(-8)		&	4.2(+17)		&	23	&	2.1(-6)		&	7.5(+16)	\\
5	&	51	&	6.4(-11)	&	1.1(+18)		&	45	&	1.6(-10)	&	1.8(+19)		&	27	&	1.3(-7)		&	4.0(+17)	\\
6	&	94	&	5.5(-12)	&	2.2(+18)		&	67	&	3.2(-10)	&	7.5(+18)		&	43	&	1.3(-9)		&	8.9(+17)	\\
7	&	110	&	4.0(-12)	&	5.9(+18)		&	73	&	1.3(-11) 	&	2.1(+18)		&	69	&	1.1(-10)	&	9.1(+18)	\\
8	&	146	&	1.8(-12) 	&	9.7(+18)		&	95	&	1.1(-11) 	&	1.5(+19)		&	70	&	2.3(-11)	&	1.6(+19)	\\
9	&	153	&	2.3(-13)	&	3.3(+18)		&	120	&	1.5(-12)	&	1.4(+19)		&	77	&	6.4(-12)	&	4.9(+18)	\\
10	&	261	&	9.4(-14)	&	1.4(+19)		&	127	&	5.8(-13)	&	7.6(+18)		&	132	&	1.8(-12) 	&	2.2(+18)	\\
11	&		&				&					&	203	&	3.3(-13)	&	3.1(+19)		&	219	&	2.6(-12)	&	1.3(+19)	\\
12	&		&				&					&	239	&	1.7(-13)	&	6.7(+19)		&	281	&	1.4(-13)	&	3.2(+19)	\\
13	&		&				&					&	254	&	6.0(-14)	&	1.4(+19)		&		&				&				\\
\hline
\end{tabular}
\end{footnotesize}
\end{table}


\subsection{Example 2: Pantograph DDE}

Consider the following pantograph differential equation (see also \cite{Brunner01}).

\begin{equation}
y'(x)= -y(x) + \frac{q}{2}y(qx) - \frac{q}{2}e^{-qx},  \qquad y(0)=1, \qquad 0 \leq x \leq T, \qquad 0<q<1.
\label{DDE_prob2}
\end{equation}

whose solution is $y_{EX}(x)= e^{-x}$. \\

Numerical methods for DDEs like (\ref{DDE_prob2}) are a topical subject of research because of two features associated to a proportional delay of the form $\tau(x)=(1-q)x$,  $0<q<1$, namely: it vanishes at $x=0$ and becomes unbounded as $x \shortrightarrow \infty$. The former one leads to difficulties in carrying out the integration of the first step, while the latter entails the need for a vast amount of computer memory if long term integration ($T>>0$) is required. In what follows we set $T=10$. 
Table \ref{test_prob2} compares the MQCM with DDESD and with the specialized reference method (REF) 
in \cite{Brunner01}, which works on a specific (\textit{geometric}) kind of mesh in order to attain superconvergence. 
The listed results for DDESD are not the best attainable, but those for which the CPU time is comparable to that of the MQCM.

\begin{table}[!hbp]	
\centering
\caption{Comparison to other methods (Example 2)}
\label{test_prob2}   
\begin{footnotesize}
\begin{tabular}{|l|ll|lll|lll|}
\hline
& \multicolumn{2}{|c|}{REF} & \multicolumn{3}{|c|}{MQ} & \multicolumn{3}{|c|}{DDESD} \\
\hline
$q$ & $DoF$ & $MAX(\epsilon)$ & $DoF$ & $MAX(\epsilon)$ & $CPU$ & $DoF$ & $MAX(\epsilon)$ & $CPU$\\
\hline
0.9 & 	1600	&	8.8(-13)	&	179	&	1.7(-13) &	9.7	&	4136	&	3.7(-14)	&	11.8	\\ 
0.5 & 	1600	&	1.2(-11)	&	135	&	2.8(-13) &	5.0	&	4136	&	3.5(-14)	&	11.2	\\ 
0.2 & 	1600	&	1.5(-11)	&	192	&	2.0(-13) & 	9.3	&	4136	&	3.3(-14)	&	11.0	\\ 
\hline
\end{tabular}
\end{footnotesize}
\end{table}

A recent improvement to the reference method \cite{Brunner01} is \cite{Ishiwataa07}
, more efficient than the former in case that long integration times $T$ are required. 
For (\ref{DDE_prob2}) with $T=10$ and $q=0.5$, it attains $|y(x=T)-y_{EX}(x=T)|= 1.8(-13)$ with $1280$ nodes. The results of the MQCM throughout the first $11$ iterations of the RSA are shown in Table \ref{RSA_prob2}. 
A fifth column has been added that lists the errors (in absolute value) of the MQCM solution at $x=T$.  

\begin{table}[!hbp]	
\centering
\caption{RSA iterations (Example 2)}
\label{RSA_prob2}   
\begin{footnotesize}
\begin{tabular}{|l|llll|}
\hline
it & DoF & RMS($\epsilon$) & Condition $\sharp$ & $\epsilon(x=T)$ \\
\hline
0 &	7 &	0.01 & 1.9(+11) &	0.02 \\
1 &	12 & 2.2(-5) & 3.3(+15) &	5.3(-6) \\
2 &	14 & 2.6(-5) & 1.7(+14) &	1.3(-4) \\
3 &	16 & 4.3(-7) & 1.1(+14) &	7.3(-7) \\
4 &	23 & 1.4(-8) & 4.2(+16) &	3.5(-9) \\
5 &	30 & 6.0(-10) & 2.2(+17) &	3.0(-10) \\
6 &	49 & 2.3(-10) & 1.1(+18) &	1.2(-10) \\
7 &	64 & 2.6(-11) & 7.3(+18) &	1.8(-10) \\
8 &	65 & 3.2(-11) & 3.3(+18) &	3.5(-11) \\
9 &	77 & 1.7(-12) & 1.2(+18) &	2.4(-12) \\
10 & 135 & 1.3(-13) & 3.9(+18) & 2.8(-13) \\
11 & 177 & 2.2(-13)  & 8.5(+18) & 8.7(-15) \\
\hline
CPU & 7.8 & & & \\
\hline
\end{tabular}
\end{footnotesize}
\end{table}

\subsection{Example 3: DDE with discontinuity propagation}
In the event that the solution $y(x)$ to the DDE has discontinuities or low-order derivative 
singularities, Kansa's method performs relatively poorly. The reason is that nonsmooth features do not belong to 
the interpolation space, which is made up of infinitely derivable MQs that cannot possibly capture 
them accurately. Any attempt to do so will bring about Gibbs' oscillations around the singularities, 
whose amplitude will not be damped by letting $N \rightarrow \infty$. An interesting approach is to 
include MQs with $c=0$ close to the 
singularities as in \cite{Jung07}. However, although the oscillations are indeed reduced, we have not been 
able to recover the high convergence rate attained with smooth solutions. In order to solve DDEs 
with piecewise smooth solutions, Kansa's method can still be applied sequentially if the domain is 
partitioned into subintervals which are $C^{\infty}$. For instance, assume that it can be predicted that 
the only three singularities take 
place at $a<x_1<x_2<x_3<b$. First, the DDE is solved in the subdomain $a \le x < x_1$ to yield 
$y_{app}^{(1)}(x)$. Then, $y_{app}^{(1)}(x)$ is used as history function for the second subdomain 
$x_1 < x < x_2$ yielding $y_{app}^{(2)}(x)$, and so on. \\

The next example (\cite{Ito91}, Example 4, also \cite{Paul92} 1.1.12), deals with a DDE having piecewise $C^{\infty}$ initial function:

\begin{equation}
y'(x)=y(x)+y(x-1)
\end{equation}

\begin{equation}
y(x)= \left\{
\begin{array}{ll}
0,	&	x \in [-1,-1/3) \nonumber\\
1, 	&	x \in [-1/3,0] \nonumber\\
\end{array}
\right.
\end{equation}

The analytical solution for $x \in [0,8/3]$ is given by
\begin{equation}
y_{EX}(x)= \left\{
\begin{array}{lr}
e^x, 						&			x \in [0,2/3], \nonumber\\
-1+C_1e^x,					&			x \in [2/3,1], \nonumber\\
xe^{x-1}+C_2e^x,			&			x \in [1,5/3], \nonumber\\
1+C_1xe^{x-1}+C_3e^x,		&			x \in [5/3,2], \nonumber\\
(\frac{x^2}{2}-x)+C_2xe^{x-1}+C_4e^x &	x \in [2,8/3]
\end{array}
\right.
\label{F:solEx3}
\end{equation}

where $C_1=1+e^{-2/3}$,$C_2=-2e^{-1}+C_1$,$c_3=\frac{5}{3}e^{-1}+C_2-e^{-5/3}-C_1\frac{5}{3}e^{-1}$, 
and $C_4= e^{-2}+2C_1e^{-1}+C_3-2C_2e^{-1}$. \\
 
Notice that the discontinuity of the initial function propagates in $x$, giving rise to singularities
of order $k$ at points $x_k= -1/3 + k,\,k\geq 0$. For both the MQCM and the SPC to cope with this problem, the integration domain $[a,b]= [0,8/3]$ must be divided into the 5 smooth subintervals in (\ref{F:solEx3}).  
The results of MQCM, SPC and DDE23 are listed in Table (\ref{test_prob3}).  
DDE23 has a \textit{'Jumps'} option which has been set to a vector that contains the locations of the discontinuities.

\begin{table}[!hbp]	
\centering
\caption{Comparison to other methods (Example 3)}
\label{test_prob3}   
\smallskip
\begin{footnotesize}
\begin{tabular}{|l|lll|}
\hline
x & $\epsilon_{MQ}$ & $\epsilon_{DDE23}$ & $\epsilon_{SPC}$ \\
\hline
0.25	&	2.0(-14) 	&	6.0(-15)	&	9.2(-13)	\\
0.5		&	4.4(-14)	&	2.2(-14)	&	1.0(-13)	\\
0.75	&	4.0(-14)	&	1.9(-13)	&	5.6(-15)	\\
1		&	6.2(-14)	&	2.6(-13)	&	1.3(-15)	\\
1.25	&	1.0(-13)  	&	3.6(-13)	&	1.1(-11)	\\
1.5		&	7.7(-14) 	&	4.9(-13)	&	1.2(-11)	\\
1.75	&	1.3(-13) 	&	6.3(-13)	&	5.8(-14)	\\
2		&	1.8(-13) 	&	8.4(-13)	&	3.3(-15)	\\
2.25	&	4.8(-13) 	&	1.5(-12)	&	7.3(-11)	\\
2.5		&	6.6(-13)  	&	2.4(-12)	&	7.9(-11)	\\
\hline
DoF		&	342		&	38386 &	 45  \\
\hline
RMS($\epsilon$) &  3.2(-13) &	9.3(-13) & \\
\hline
CPU &  11 &	246 &  \\
\hline
\end{tabular}
\end{footnotesize}
\end{table}


\subsection{Example 4: Neutral state-delay DDE} 

Neutral DDEs (which involve lagged derivatives) are considered tougher to handle with numerical methods than
retarded ODEs and are an active research field. In the following problem,
taken from \cite{Norkin73} (see also \cite{Paul92} 2.3.4), the delay is a function of the solution itself, and
therefore the DDE is nonlinear. The MQCM tackles it with Powell's method,
implemented by the option \textit{'dogleg'} of MATLAB nonlinear solver \textit{fsolve}.  

\begin{equation}
\left\{
\begin{array}{ll}
y'(x)= -y'(y(x)-2),  \qquad x \geq 0 \\
y(x)=1-x, \qquad x\le 0.
\end{array}
\right.
\end{equation}

The exact solution is $y_{EX}(x)=1+x,\, \, 0 \le x \le 1$.

\begin{table}[!hbp]	
\centering
\caption{RSA iterations (Example 4)}
\label{T:prob4RSA}   
\smallskip
\begin{footnotesize}
\begin{tabular}{|l|llll|}
\hline
It & DoF & RMS($\epsilon$) & Condition & NL iter\\
\hline 
0 	&	7 	&	0.0452		&	7.5(+11)	&	\textit{f}\\ 
1  	&	13	&  	1.3(-11)	&	5.1(+13)	&	22\\
2  	&	14	&  	1.2(-11)	&	2.9(+15)	&	16\\
3  	&	18	&  	1.9(-13)	&	1.1(+14)	&	19\\
4  	&	24 	& 	2.0(-14)	&	2.5(+15)	&	15\\
\hline
CPU & \multicolumn{4}{l|}{7.5} \\
\hline
\end{tabular}
\end{footnotesize}
\end{table}

The initial guess of $y(x)$ required to trigger Powell's method 
is $y_{GUESS}(x)=0$.
The entry \textit{f} in Table \ref{T:prob4RSA} means that Powell's method fails
to converge in the maximum number of iterations allowed (set to 30). Nevertheless, it
yields an approximation good enough to be used as a guess for the nonlinear solution
with 13 MQs (whose solution is in turn used as a guess for the next RSA iteration, and so on). 
For reference, DDENSD yields $RMS(\epsilon_{REF})= 2.2(-9)$ in $50.1$ s. CPU time.

\subsection{Example 5: Vanishing state-delay DDE}
This example is a nonlinear neutral differential equation with vanishing state delay. It was first proposed in \cite{Enright98} as a modification of a problem originally considered in \cite{Castleton73}:

\begin{equation}
\left\{
\begin{array}{ll}
y'(x)= \cos(x)\big[ 1 + y\big(xy^2(x)\big) \big] + cy(x)y'\big( xy^2(x) \big) + g(x),    \qquad 0 \leq x \leq \pi \\
g(x)= (1-c)\sin(x)\cos\big(x\sin^2(x)\big) - \sin\big(x+x\sin^2(x)\big) \\
y(0)= 0 \\
\end{array}
\right.
\label{F:DDE5}
\end{equation}

For every choice of the parameter $c$, the exact solution is $y_{EX}(x)=sin(x)$. 
Because the delay vanishes at $x=0,\pi/2,3\pi/2,...$, the numerical solution of (\ref{F:DDE5}) by Runge-Kutta methods causes some difficulties. For the MQCM, the main difficulty is that the condition number must be kept low enough (below $10^{14}$) for the nonlinear solver to converge in a reasonable number of nonlinear iterations (again Powell's algorithm in the \textit{fsolve} routine). Therefore, we have set $N^{(0)}=11$ and $\mu=\sqrt{20/N^{(0)}}$. The initial hint of the solution is $y_{GUESS}=1/2$. As reference results ($REF$), we have taken those of \cite{Guglielmi00} (example 2), where (\ref{F:DDE5}) is solved by the Radau-type code RADAR5 (Table \ref{T:prob5comp}). In the case $c=1$, there is a singularity at $x=\pi/2$ -in the sense that $y'(\pi/2)$ is not well defined- and the MQCM with default parameters fails.   

\begin{table}[!hbp]	
\centering
\caption{RSA iterations (Example 5)}
\label{T:prob5comp}   
\smallskip
\begin{footnotesize}
\begin{tabular}{|l|llllll|}
\hline
c & DoF & Condition & CPU & RMS($\epsilon$) &	$\epsilon_{MQ}(x=\pi)$ &	$\epsilon_{REF}(x=\pi)$ \\
\hline
-1.0	&	65	&	1.8(13)	&	58.3	&	4.7(-9)	&	3.3(-9)		&	1.8(-9)	\\
-0.7	&	44	&	2.5(11)	&	20.4	&	3.2(-8)	&	9.5(-9)		&	4.2(-9)	\\
-0.3	&	44	&	9.1(10)	&	9.8		&	3.2(-8)	&	7.5(-8)		&	1.7(-10)\\
 0.0	&	69	&	6.8(12)	&	30.0	&	3.0(-8)	&	1.6(-8)		&	1.2(-9)	\\
 0.3	&	46	&	9.9(10)	&	17.1	&	4.3(-9)	&	2.5(-9)		&	1.0(-9)	\\
 0.7	&	49	&	3.0(11)	&	36.2	&	1.1(-9)	&	7.2(-10)	&	5.3(-9)	\\
 1.0	&	$\textit{f}$	&			&			&		&		&	4.3(-8)	\\
\hline
\end{tabular}
\end{footnotesize}
\end{table}

\subsection{Example 6: Second order DDE}

The last example illustrates the ability of the MQCM to accurately solve higher-order DDEs. Since
equations of this type are less common in the literature, most solvers are not
designed to handle them. In order to compare, we have transformed a system of two state-delay
DDEs into a second-order DDE:
\begin{equation}
\left\{
\begin{array}{ll}
y''(x)= \big( \exp[1-y(x)]-x\big)y\big(x-\exp [1-y(x)]\big)y'(x)^2  \qquad\qquad x\ge 1\\
y(x)= \log(x)  \qquad\qquad 0<x\le 1
\end{array}
\right.
\label{F:DDE6a}
\end{equation} 

which is obtained by derivation of $y_2(x)$ and insertion into $y_1''(x)$ in   

\begin{equation}
\left\{
\begin{array}{ll}
y_1'(x)= y_2(x) \qquad\qquad x\ge 1 \\
y_2'(x)= \big(\exp[1-y_1(x)]-x\big)y_2(x)(x-\exp [1-y_1(x)])y_2^2(x) \qquad\qquad x\ge 1\\
y_1(x)= \log(x)  \qquad\qquad 0<x\le 1 \\
y_2(x)= 1/x  \qquad\qquad 0<x\le 1
\end{array}
\right.
\label{F:DDE6b}
\end{equation}

(see \cite{AlMutib87} and \cite{Paul92}, 1.4.17). The exact solution is $y_{EX}(x)=y_{1,EX}(x)=\log(x)$.

We consider the interval $[a,b]=[1,5]$. In this problem, the second derivative of the
multiquadric (\ref{F:MQ}) is required, as well as \textit{two} extra MQ centers for
PDEBC (since $y$,$y'$, and $y''$ are enforced at $x=a$). Such
centers are placed at $x_0=a-\Delta$ and $x_{-1}=a-2\Delta$. Results are shown in Table \ref{T:prob6RSA}. 

\begin{table}[!hbp]	
\centering
\caption{RSA iterations (Example 6)}
\label{T:prob6RSA}   
\smallskip
\begin{footnotesize}
\begin{tabular}{|l|llll|}
\hline
It & DoF & RMS($\epsilon$) & Condition & NL iter\\
\hline 
0 	&	12	&	0.103	&	2.4(+10)	&	\textit{f}\\
1  	&	18	&  2.9(-5)	&	2.0(+10)	&	\textit{f}\\
2  	&	21	&  3.1(-7)	&	6.8(+11)	&	1 \\
3  	&	28	&  3.3(-8)	&	8.7(+10)	&	1 \\
4  	&	37 	&  3.7(-9)	&	5.5(+13)	&	1 \\
5  	&	49 	&  2.1(-11)	&	2.5(+13)	&	1 \\
6  	&	63 	&  6.2(-11)	&	1.7(+15)	&	1 \\
7  	&	72 	&  8.6(-12)	&	3.5(+14)	&	1 \\
\hline
CPU & \multicolumn{4}{l|}{48.5} \\
\hline
\end{tabular}
\end{footnotesize}
\end{table}

An indirect reference is provided by DDESD which solves
(\ref{F:DDE6b}) with RMS($\epsilon_{REF})= 2.3(-13)$ (for $y_1(x)$) in 2.1 s.
Notice that ill-condition must be kept lower in order to ensure
convergence of the nonlinear solver, thus limiting accuracy. A
possible alternative would be to use a Newton-type routine with
the analytical Jacobian to the given DDE.
The good performance shown by MQCM relies on the accuracy with which numerical derivatives
are reproduced in Kansa's method. While not every system of $m$ DDEs can be transformed into
a single DDM of order $m$, there are many cases where this transformation can be
an advantageous alternative for the solution of DDE systems with the MQCM.

\section{Conclusions}
A novel numerical method for the solution of DDEs has been presented. It relies on
the multiquadric collocation method introduced by Kansa combined with
the RSA algorithm by Driscoll and Heryudono for node adaptivity, which
uses residual as refinement criterium. As long
as the solution of the DDE is smooth (or piecewise smooth, with the position
of the singularities being known in advance), the present method can accurately handle
a large variety of such problems, including state-delay, neutral, and high-order
DDEs. Moreover, the scheme is straightforward to code and enjoys spectral convergence.
Because in this paper the stress is placed on simplicity, nonlinearities
are fed to a general-purpose solver, without attempting to optimize. Possible
improvements include: providing an analytical Jacobian, and linearizing
the DDE along the lines of \cite{Fasshauer01}\cite{Bernal08}.

\end{document}